\documentclass[11pt]{article}

\usepackage[letterpaper,margin=1.08in,headheight=14pt]{geometry}
\usepackage[T1]{fontenc}
\usepackage[utf8]{inputenc}
\usepackage{lmodern}
\usepackage{microtype}
\usepackage{amsmath,amssymb,amsthm,mathtools}
\usepackage[dvipsnames]{xcolor}
\usepackage{fancyhdr}
\usepackage{xurl}
\usepackage{hyperref}
\usepackage{aliascnt}
\usepackage{comment}
\usepackage[nameinlink,capitalize,noabbrev]{cleveref}

\definecolor{deepblue}{HTML}{173B57}
\definecolor{mutedblue}{HTML}{476C87}
\hypersetup{
  colorlinks=true,
  linkcolor=deepblue,
  citecolor=mutedblue,
  urlcolor=deepblue,
  pdfauthor={Weibo Fu, Yanjun Han, Guanyang Wang, Jun Yan, Peng Zhang, and Zhengqing Zhou},
  pdftitle={Sharp small-deviation inequalities for sums of independent nonnegative random variables},
  pdfsubject={Sharp small-deviation inequalities for sums of independent nonnegative random variables},
  pdfkeywords={Feige conjecture, Gaffke conjecture, Dirichlet calibration, convex geometry, Grunbaum inequality}
}

\allowdisplaybreaks

\newtheorem{theorem}{Theorem}[section]
\newaliascnt{lemma}{theorem}
\newtheorem{lemma}[lemma]{Lemma}
\aliascntresetthe{lemma}
\newaliascnt{corollary}{theorem}
\newtheorem{corollary}[corollary]{Corollary}
\aliascntresetthe{corollary}
\theoremstyle{remark}
\newaliascnt{remark}{theorem}

\aliascntresetthe{remark}
\crefname{lemma}{Lemma}{Lemmas}
\crefname{corollary}{Corollary}{Corollaries}
\crefname{remark}{Remark}{Remarks}

\numberwithin{equation}{section}

\newcommand{\E}{\mathbb E}
\newcommand{\Pp}{\mathbb P}
\newcommand{\R}{\mathbb R}
\newcommand{\vol}{\operatorname{vol}}
\newcommand{\Dir}{\operatorname{Dir}}

\title{Sharp small-deviation inequalities for sums of independent\\nonnegative random variables}
\author{
\begin{tabular}{ccc}
\shortstack{Weibo Fu\\\texttt{wfu@math.princeton.edu}}
&
\shortstack{Yanjun Han\\\texttt{yanjunhan@nyu.edu}}
&
\shortstack{Guanyang Wang\\\texttt{guanyang.wang@rutgers.edu}}
\\[1.2em]
\shortstack{Jun Yan\\\texttt{junyan65@stanford.edu}}
&
\shortstack{Peng Zhang\\\texttt{pz149@rutgers.edu}}
&
\shortstack{Zhengqing Zhou\\\texttt{zqzhou@stanford.edu}}
\end{tabular}
}
\date{26 July 2026}

\begin{document}

\maketitle
\thispagestyle{plain}

\begin{abstract}
Let $(X_1,\ldots,X_n)$ be independent nonnegative random variables with $\E X_i\le1$, and write $S=\sum_iX_i$. For $\delta>0$, we prove that
\[
 \Pp\left(S<\E S+\delta\right)\ge b_{n,\delta} :=
 \begin{cases}
 \delta(\frac{n}{n+\delta})^n,&0<\delta<1,\\
 (1-\frac1{n+\delta})^n,&\delta\ge1.
 \end{cases}
\]
The bound is sharp for every $n$ and $\delta\ge 1$. In particular, since $b_{n,\delta} \ge e^{-1}$ for $\delta \ge 1$, our result proves Feige's conjecture \cite{Feige2004} in the affirmative for $\delta\ge 1$.

The proof is found by ChatGPT 5.6 Pro. It combines the exact Dirichlet calibration theorem of Vlassis and Thomas \cite{VlassisThomas2026}, which resolves Gaffke's conjecture in statistics, with results in convex geometry including Gr\"unbaum's centroid theorem \cite{Grunbaum1960} and its generalization by Letwin and Yaskin \cite{LetwinYaskin2024}.
\end{abstract}

\section{Introduction}\label{sec:introduction}

We begin by stating our main result, which gives a sharp lower bound for sums of independent nonnegative random variables.

\begin{theorem}[Main theorem]\label{thm:main}
Let $n\ge1$, let $X_1,\ldots,X_n$ be independent nonnegative random variables with $\E X_i\le1$, and put $S=\sum_iX_i$. For $\delta>0$, set
\begin{equation}\label{eq:beta}
 b_{n,\delta}=
 \begin{cases}
 \displaystyle \delta\left(\frac{n}{n+\delta}\right)^n,&0<\delta<1,\\[0.8em]
 \displaystyle \left(1-\frac1{n+\delta}\right)^n,&\delta\ge1.
 \end{cases}
\end{equation}
Then
\begin{equation}\label{eq:main}
 \Pp\left(S<\E S+\delta\right)\ge b_{n,\delta}.
\end{equation}
For every $n$ and every $\delta\ge1$, the bound is optimal.
\end{theorem}

The following corollary is immediate. 
\begin{corollary}[Universal lower bound]\label{cor:main}
Under the conditions of \Cref{thm:main}, it holds that
\[
\Pp\left(S<\mathbb{E}S+\delta\right)>
\begin{cases}
\delta e^{-\delta}, & 0<\delta<1,\\[2mm]
e^{-1}, & \delta\ge 1.
\end{cases}
\]
\end{corollary}

Feige \cite{Feige2004} proposed two closely related conjectures. The general form conjectured that, for every $\delta>0$, the sharp lower bound is \cite{Feige2004}
\[
 \min\left\{\frac{\delta}{1+\delta},e^{-1}\right\}.
\]
The more commonly stated version, corresponding to $\delta=1$, asserts the optimal universal bound
\[
 \Pp\left(S<\E S+1\right)\ge e^{-1}.
\]
\Cref{thm:main} and \Cref{cor:main} prove the conjectured sharp bound throughout the range $\delta\ge1$, and obtain a weaker constant $\delta e^{-\delta}$ compared to Feige's conjecture $\min\{\frac{\delta}{1+\delta},e^{-1}\}$ for $0<\delta<1$. 

\subsection{Background and related work}

The problem is a natural sharp small-deviation question for sums of independent nonnegative random variables under only first-moment constraints. Its study goes back to Samuels' work on Markov-type inequalities \cite{Samuels1966,Samuels1968,Samuels1969}. Feige isolated the unit-slack case and proved the first universal lower bound, $1/13$ \cite{Feige2004}. This was later improved to $1/8$ by He, Zhang, and Zhang \cite{HeZhangZhang2010}, to approximately $0.14$ by Garnett \cite{Garnett2020}, and to $0.1798$ by Guo, He, Ling, and Liu \cite{GuoHeLingLiu2020}. Related formulations and inequalities were studied by Elton \cite{Elton2009}, Oleszkiewicz \cite{Oleszkiewicz2012}, and Paulin \cite{Paulin2017}. Feige's $e^{-1}$ conjecture was previously proved under additional assumptions: for discrete log-concave distributions by Alqasem, Aravinda, Marsiglietti, and Melbourne \cite{AlqasemEtAl2024}, and for identically distributed random variables by Egozcue and Fuentes Garc\'ia \cite{EgozcueFuentes2025}. 

\subsection{Proof Outline}

\paragraph{Sharpness for $\delta\ge1$.}
Let the $X_i$ be independent and satisfy
\[
 X_i=\begin{cases}
 n+\delta,&\text{with probability }1/(n+\delta),\\
 0,&\text{otherwise}.
 \end{cases}
\]
Then $\E X_i=1$, and $S<\E S+\delta=n+\delta$ holds exactly when all the $X_i$ vanish. Hence
\[
 \Pp\left(S<\E S+\delta\right)=\left(1-\frac1{n+\delta}\right)^n.
\]
For $\delta\ge1$, this is the second branch of \eqref{eq:beta}.

\paragraph{Proof sketch of \cref{thm:main}.}
The proof combines the recent solution of Gaffke's conjecture by Vlassis and Thomas \cite{VlassisThomas2026} with Gr\"unbaum's centroid theorem \cite{Grunbaum1960}. Gaffke introduced the statistic below in 2005 as a finite-sample, distribution-free test for a one-sided mean hypothesis and conjectured its validity under independent coordinatewise mean constraints \cite{Gaffke2005}. He established the corresponding i.i.d. asymptotics, reduced the finite-sample problem at each level to mean-one two-point laws, proved the case $n=2$, and reported numerical verification through $n=15$ \cite{Gaffke2005,VlassisThomas2026}. The confidence interval obtained by inverting the test for bounded observations was later studied by Learned-Miller and Thomas \cite{LearnedMillerThomas2020}. The auxiliary chain measure for the two-point argument used in the proof of \cite[Proposition~3]{VlassisThomas2026} is precisely a chain distribution underlying the Lov{\'a}sz extension \cite{Lovasz1983,Dughmi2009}. After Vlassis and Thomas proved the conjecture in full generality, Ming, Ramdas, Shen, Wang, and Waudby-Smith showed that the resulting $p$-value is inadmissible for $n\ge2$, although the corresponding equal-tail confidence interval is first-order asymptotically efficient \cite{MingEtAlGaffke2026}.

Let $D=(D_0,\ldots,D_n)$ be uniformly distributed on the standard $n$-simplex $\Delta_n$, and define
\begin{equation}\label{eq:K-intro}
 K_n(x)=\Pp_D\left(\sum_{i=1}^n x_iD_i\le1\right),
 \qquad x\in[0,\infty)^n.
\end{equation}
The theorem of Vlassis and Thomas states that, if $X=(X_1,\ldots,X_n)$ has independent nonnegative coordinates satisfying $\E X_i\le1$ for every $i$, then
\[
 \Pp\left(K_n(X)\le\alpha\right)\le\alpha,
 \qquad 0\le\alpha\le1.
\]

The geometric part upper bounds the function $K_n(x)$, which is the normalized $n$-volume of the halfspace $\{D: \sum_{i=1}^n x_iD_i \le 1\}\cap \Delta_n$. Let $g=(g_0,\dots,g_n)=(\frac{1}{n+1},\dots,\frac{1}{n+1})$ be the centroid of the $n$-simplex, then if $\sum_i x_i\ge n+\delta$, 
\begin{align*}
K_n(x) \le \Pp_D\left(\sum_{i=1}^n x_iD_i\le \frac{n+1}{n+\delta}\sum_{i=1}^n x_ig_i \right). 
\end{align*}
When $\delta = 1$, the hyperplane $\{D: \sum_{i=1}^n x_iD_i = 1\}$ contains the centroid $g$, so that Gr\"unbaum inequality \cite{Grunbaum1960} shows that $K_n(x)\le 1-(\frac{n}{n+1})^n\le 1-e^{-1}$. For general $\delta>0$, the generalized Gr\"unbaum inequality of Letwin and Yaskin \cite[Theorem~4]{LetwinYaskin2024} gives
\[
K_n(x) \le 1-b_{n,\delta}. 
\]
Combining the above geometric estimates with the Vlassis--Thomas theorem gives \eqref{eq:main}.

\section*{Statement on AI use}

The initial proof is found by ChatGPT 5.6 Pro. The authors subsequently checked, revised, and rewrote the argument, and take full responsibility for the final content.

An accompanying Lean formalization, developed with Codex and available at \url{https://github.com/pengzhang91/Feige}, provides an end-to-end formal proof of Feige's $e^{-1}$ conjecture. It formalizes the main theorem of Vlassis and Thomas \cite{VlassisThomas2026}, Gr\"unbaum's centroid theorem \cite{Grunbaum1960}, and the final argument that combines these two results.

\section{Proof of the main theorem}\label{sec:proof}

\subsection{Preliminaries}
We state several preliminary results that will be useful ingredients in our proof. The first ingredient is a calibration result shown recently by Vlassis and Thomas \cite{VlassisThomas2026}. Let
\[
 \Delta_n=\left\{d=(d_0,\ldots,d_n)\in[0,\infty)^{n+1}:
 \sum_{i=0}^n d_i=1\right\}
\]
be the probability $n$-simplex, and let $D\sim\Dir(1,\ldots,1)$, the uniform probability law on $\Delta_n$. For $y\in[0,\infty)^n$, set
\begin{equation}\label{eq:K}
 K_n(y)=\Pp_D\left(\sum_{i=1}^n y_iD_i \le 1\right).
\end{equation}
We use the following theorem as an external input \cite[Theorem~1]{VlassisThomas2026}.

\begin{theorem}[Vlassis--Thomas]\label{thm:VT}
If $Y_1,\ldots,Y_n$ are independent nonnegative random variables with $\E Y_i\le1$, then
\begin{equation}\label{eq:VT}
 \Pp\left(K_n(Y_1,\ldots,Y_n)\le\alpha\right)\le\alpha,
 \qquad 0\le\alpha\le1.
\end{equation}
\end{theorem}

The next ingredient is a recent generalization of Gr\"unbaum inequality by Letwin and Yaskin \cite{LetwinYaskin2024}. 

\begin{theorem}[Letwin--Yaskin]\label{thm:letwin-yaskin}
Let $K\subseteq \mathbb{R}^n$ be a convex body with centroid at the origin. Then for $\alpha\in (-1,\frac{1}{n})$, every non-zero vector $\xi\in \mathbb{R}^n$, and the closed halfspace
\begin{align*}
H_{\alpha,\xi} := \left\{x\in \mathbb{R}^n: \langle x,\xi \rangle \ge \alpha h_K(-\xi)\right\}, 
\end{align*}
where $h_K(t)=\sup_{x\in K}\langle x,t\rangle$ is the support function of $K$, it holds that
\begin{align*}
\frac{\vol_n(K\cap H_{\alpha,\xi})}{\vol_n(K)} \ge \begin{cases}
    (\frac{n-\alpha}{n+1})^n &\text{if }-1<\alpha\le 0, \\
    (\frac{n}{n+1})^n(\alpha+1)^{n-1}(1-\alpha n) &\text{if } 0<\alpha<\frac{1}{n}.
\end{cases}
\end{align*}
\end{theorem}

When $\alpha=0$, the lower bound $(\frac{n}{n+1})^n$ is the celebrated inequality of Gr\"unbaum \cite{Grunbaum1960} in convex geometry. 

\subsection{Proof of \Cref{thm:main}}
We first establish an upper bound of $K_n(y)$ in \eqref{eq:K} when $\sum_{i=1}^n y_i \ge n+\delta$. To this end, we apply \Cref{thm:letwin-yaskin} to the shifted $n$-simplex $K = \Delta_n - (\frac{1}{n+1},\dots,\frac{1}{n+1})\subseteq \mathbb{R}^{n+1}$, with centroid at the origin. Since $K\subseteq {\bf 1}^\perp$, we may identify $K$ as a convex body in $\mathbb{R}^n$; under this identification, every direction $\xi \in \mathbb{R}^{n+1}$ can be converted into a direction in $\mathbb{R}^n$ by the orthogonal projection of $\xi$ to ${\bf 1}^\perp$, which does not change either $\langle x,\xi\rangle$ for $x\in K$ or $h_K(-\xi)$ in \Cref{thm:letwin-yaskin}. 

We choose the direction $\xi = (0,y_1,\dots,y_n)$. Since $y_i\ge 0$, 
\begin{align*}
h_K(-\xi) = \max_{x\in \Delta_n} \sum_{i=1}^n \left(\frac{1}{n+1}-x_i\right)y_i = \frac{1}{n+1}\sum_{i=1}^n y_i, 
\end{align*}
with maximum attained at $x=(1,0,\dots,0)$. Therefore, if  $\sum_{i=1}^n y_i \ge n+\delta$, then
\begin{align*}
1-K_n(y) = \Pp_D\left(\sum_{i=1}^n y_iD_i > 1\right) &=\Pp_D\left(\sum_{i=1}^n y_iD_i\ge 1\right) \\
&\ge \Pp_D\left(\sum_{i=1}^n y_iD_i\ge \frac{1}{n+\delta}\sum_{i=1}^n y_i \right)\\
&= \Pp_D\left(\sum_{i=0}^n \xi_i\left(D_i-\frac{1}{n+1}\right)\ge \frac{1-\delta}{n+\delta} h_K(-\xi)\right) \\
&= \frac{\vol_n(K\cap H_{(1-\delta)/(n+\delta),\xi})}{\vol_n(K)}\\
&\ge \begin{cases}
    \delta(\frac{n}{n+\delta})^n &\text{if }0<\delta < 1 \\
    (1-\frac{1}{n+\delta})^n &\text{if } \delta \ge 1
\end{cases}. 
\end{align*}
Here the second equality follows because $\left\{d\in\Delta_n:\sum_{i=1}^n y_i d_i=1\right\}
$ has zero $n$-dimensional simplex volume, and the final inequality follows from \Cref{thm:letwin-yaskin} after simple algebra. Therefore, we have shown that $K_n(y)\le 1-b_{n,\delta}$ whenever $\sum_{i=1}^n y_i \ge n+\delta$. 

Next we prove \Cref{thm:main}. Let $\mu_i=\E X_i$ and define
\[
 Y_i=X_i+1-\mu_i.
\]
The variables $Y_i$ are independent and nonnegative, $\E Y_i=1$, and
\[
 \sum_{i=1}^nY_i=S+n-\E S.
\]
The event $\{S\ge\E S+\delta\}$ implies $\sum_iY_i\ge n+\delta$, and further $K_n(Y)\le 1-b_{n,\delta}$. 
Therefore, by \cref{thm:VT},
\begin{align*}
 \Pp\left(S\ge\E S+\delta\right)
 \le\Pp\left(K_n(Y)\le1-b_{n,\delta}\right)
 \le1-b_{n,\delta}.
\end{align*}
Taking complements proves \eqref{eq:main}.

\bibliographystyle{alpha}
\small
\bibliography{references}

@article{AlqasemEtAl2024,
  author  = {Alqasem, Abdulmajeed and Aravinda, Heshan and Marsiglietti, Arnaud and Melbourne, James},
  title   = {On a conjecture of {Feige} for discrete log-concave distributions},
  journal = {SIAM Journal on Discrete Mathematics},
  volume  = {38},
  number  = {1},
  pages   = {93--102},
  year    = {2024},
  doi     = {10.1137/22M1539514}
}

@article{Dughmi2009, 
    author = {Dughmi, Shaddin}, 
    title = {Submodular functions: Extensions, distributions, and algorithms. A survey}, 
    year = {2009},
    eprint = {0912.0322},
    archivePrefix = {arXiv},
    primaryClass  = {math.PR},
    note          = {arXiv:0912.0322}
}

@article{EgozcueFuentes2025,
  author        = {Egozcue, Mart{\'i}n and Fuentes Garc{\'i}a, Luis},
  title         = {A short proof of {Feige}'s conjecture for identically distributed random variables},
  year          = {2025},
  eprint        = {2509.19949},
  archivePrefix = {arXiv},
  primaryClass  = {math.PR},
  note          = {arXiv:2509.19949}
}

@article{Elton2009,
  author        = {Elton, John H.},
  title         = {Notes on {Feige}'s gumball machines problem},
  year          = {2009},
  eprint        = {0908.3528},
  archivePrefix = {arXiv},
  primaryClass  = {math.PR},
  note          = {arXiv:0908.3528}
}

@inproceedings{Feige2004,
  author    = {Feige, Uriel},
  title     = {On sums of independent random variables with unbounded variance, and estimating the average degree in a graph},
  booktitle = {Proceedings of the Thirty-Sixth Annual ACM Symposium on Theory of Computing},
  pages     = {594--603},
  year      = {2004}
}

@article{Gaffke2005,
  author  = {Gaffke, Norbert},
  title   = {Three test statistics for a nonparametric one-sided hypothesis on the mean of a nonnegative variable},
  journal = {Mathematical Methods of Statistics},
  volume  = {14},
  number  = {4},
  pages   = {451--467},
  year    = {2005}
}

@article{Garnett2020,
  author  = {Garnett, Brian},
  title   = {Small deviations of sums of independent random variables},
  journal = {Journal of Combinatorial Theory, Series A},
  volume  = {169},
  pages   = {105119},
  year    = {2020},
  doi     = {10.1016/j.jcta.2019.105119}
}

@article{Grunbaum1960,
  author  = {Gr{\"u}nbaum, Branko},
  title   = {Partitions of mass-distributions and of convex bodies by hyperplanes},
  journal = {Pacific Journal of Mathematics},
  volume  = {10},
  number  = {4},
  pages   = {1257--1261},
  year    = {1960},
  doi     = {10.2140/pjm.1960.10.1257}
}

@article{GuoHeLingLiu2020,
  author        = {Guo, Jiayi and He, Simai and Ling, Zi and Liu, Yicheng},
  title         = {Bounding probability of small deviation on sum of independent random variables: Combination of moment approach and {Berry--Esseen} theorem},
  year          = {2020},
  eprint        = {2003.03197},
  archivePrefix = {arXiv},
  primaryClass  = {math.PR},
  note          = {arXiv:2003.03197}
}

@article{HeZhangZhang2010,
  author  = {He, Simai and Zhang, Jiawei and Zhang, Shuzhong},
  title   = {Bounding probability of small deviation: A fourth moment approach},
  journal = {Mathematics of Operations Research},
  volume  = {35},
  number  = {1},
  pages   = {208--232},
  year    = {2010},
  doi     = {10.1287/moor.1090.0438}
}

@article{LearnedMillerThomas2020,
  author        = {Learned-Miller, Erik and Thomas, Philip S.},
  title         = {A new confidence interval for the mean of a bounded random variable},
  year          = {2020},
  eprint        = {1905.06208},
  archivePrefix = {arXiv},
  primaryClass  = {stat.ME},
  note          = {arXiv:1905.06208v2}
}

@incollection{Lovasz1983,
    author = {Lov{\'a}sz, L{\'a}szl{\'o}}, 
    title = {Submodular functions and convexity}, 
    booktitle = {Mathematical Programming: The State of the Art}, 
    editor = {Bachem, Achim and Gr{\"o}tschel, Martin and Korte, Bernhard}, 
    pages = {235--257}, 
    publisher = {Springer}, 
    address = {Berlin, Heidelberg}, 
    year = {1983}, 
    doi = {10.1007/978-3-642-68874-4_10} 
}

@article{MingEtAlGaffke2026,
  author        = {Ming, Jiahao and Ramdas, Aaditya and Shen, Yi and Wang, Ruodu and Waudby-Smith, Ian},
  title         = {{Gaffke}'s confidence interval for the mean of bounded data is inadmissible but asymptotically efficient},
  year          = {2026},
  eprint        = {2607.18661},
  archivePrefix = {arXiv},
  primaryClass  = {math.ST},
  note          = {arXiv:2607.18661v1, posted 21 July 2026}
}

@incollection{Oleszkiewicz2012,
  author    = {Oleszkiewicz, Krzysztof},
  title     = {On some extension of {Feige}'s inequality},
  booktitle = {Geometric Aspects of Functional Analysis},
  editor    = {Klartag, Bo'az and Mendelson, Shahar and Milman, Vitali D.},
  series    = {Lecture Notes in Mathematics},
  volume    = {2050},
  pages     = {345--353},
  publisher = {Springer},
  address   = {Berlin},
  year      = {2012},
  doi       = {10.1007/978-3-642-29849-3_21}
}

@article{Paulin2017,
  author        = {Paulin, Roland},
  title         = {On some conjectures of {Samuels} and {Feige}},
  year          = {2017},
  eprint        = {1703.05152},
  archivePrefix = {arXiv},
  primaryClass  = {math.PR},
  note          = {arXiv:1703.05152}
}

@article{Samuels1966,
  author  = {Samuels, Stephen M.},
  title   = {On a {Chebyshev}-type inequality for sums of independent random variables},
  journal = {Annals of Mathematical Statistics},
  volume  = {37},
  number  = {1},
  pages   = {248--259},
  year    = {1966},
  doi     = {10.1214/aoms/1177699614}
}

@techreport{Samuels1968,
  author      = {Samuels, Stephen M.},
  title       = {More on a {Chebyshev}-type inequality for sums of independent random variables},
  institution = {Department of Statistics, Purdue University},
  number      = {155},
  year        = {1968}
}

@article{Samuels1969,
  author  = {Samuels, Stephen M.},
  title   = {The {Markov} inequality for sums of independent random variables},
  journal = {Annals of Mathematical Statistics},
  volume  = {40},
  number  = {6},
  pages   = {1980--1984},
  year    = {1969}
}

@article{VlassisThomas2026,
  author        = {Vlassis, Nikos and Thomas, Philip S.},
  title         = {An exact distribution-free test for means of nonnegative random variables},
  year          = {2026},
  eprint        = {2607.08415},
  archivePrefix = {arXiv},
  primaryClass  = {math.ST},
  note          = {arXiv:2607.08415v1, posted 9 July 2026}
}

@article{LetwinYaskin2024,
  author        = {Letwin, Brayden and Yaskin, Vladyslav},
  title         = {A generalization of {Gr\"unbaum}'s inequality},
  year          = {2024},
  eprint        = {2410.04741},
  archivePrefix = {arXiv},
  primaryClass  = {math.MG},
  note          = {arXiv:2410.04741; to appear in Israel Journal of Mathematics}
}

\end{document}